\documentclass[12pt,oneside]{amsart}
\usepackage{latexsym}
\usepackage{amsmath,amssymb,amsthm}
\usepackage{amscd}
\usepackage[dvips]{graphics}
\newcommand{\rp}{\mathbb{RP}}
\newcommand{\re}{\mathbb{R}}
\newcommand{\co}{\mathbb {C}}
\newcommand{\cp}{\mathbb {CP}}

\newcommand{\pgl}[1]{\mathbf{PGL}(#1,\mathbb{R})}
\newcommand{\Sl}[1]{\mathbf{SL}(#1,\mathbb{R})}
\newcommand{\sa}[1]{\mathbf{SA}(#1,\mathbb{R})}
\newcommand{\so}[1]{\mathbf{SO}(#1,\mathbb{R})}

\newcommand{\na}{\nabla}
\newcommand{\sfrac}[2]{{\textstyle \frac{#1}{#2}}}

\newtheorem{prop}{Proposition} 
\newtheorem{cor}[prop]{Corollary}
\newtheorem{thm}{Theorem}

\theoremstyle{remark}
\newtheorem*{rem}{Remark}

\begin{document}

\title{Survey on Affine Spheres}
\author{John Loftin}
\maketitle

\section{Introduction}
Affine spheres were introduced by \c{T}i\c{t}eica in
\cite{tzitzeica08,tzitzeica09}, and studied later by Blaschke,
Calabi, and Cheng-Yau, among others. These are hypersurfaces in
affine $\re^{n+1}$ which are related to real Monge-Amp\`ere
equations, to projective structures on manifolds, and to the
geometry of Calabi-Yau manifolds.  In this survey article, we will
outline the theory of affine spheres their relationships to these
topics.

Affine differential geometry is the study of those differential
properties of hypersurfaces of $\re^{n+1}$ which are invariant under
all volume-preserving affine transformations.  Affine differential
geometry is largely traced to \c{T}i\c{t}eica's papers in 1908-09,
although for curves in $\re^2$, one of the main invariants, the
affine normal, was already introduced by Transon \cite{transon41} in
1841. Given a smooth hypersurface $H\subset \re^{n+1}$, the affine
normal $\xi$ is an affine-invariant transverse vector field to $H$.
Define the special affine group as
$$\sa{n+1} = \{\Phi\!: x \mapsto Ax+b, \,\det A = 1\}.$$  The
invariance property of the affine normal is then $$\Phi_*\xi_H(x) =
\xi_{\Phi(H)} (\Phi(x))$$ for any $x\in H$.  An \emph{improper
affine sphere} is a hypersurface $H$ whose affine normals are all
parallel, while a \emph{proper affine sphere} is a hypersurface
whose affine normal lines all meet in a point, the \emph{center} of
the affine sphere.  By symmetry, a Euclidean sphere must be a proper
affine sphere, and affine invariance then shows that all ellipsoids
are affine spheres also.  More generally, quadric hypersurfaces are
the canonical examples of affine spheres.

We will mainly focus on the case of convex hypersurfaces in this
survey, since in this case the natural invariant metric, the affine
metric, is positive definite, and so we can exploit techniques of
Riemannian geometry and elliptic PDEs. We assume for a convex
hypersurface that the affine normal points to the convex side of the
hypersurface. Improper affine sphere are then called \emph{parabolic
affine spheres}, and the primary example is an elliptic paraboloid.
For convex hypersurfaces, there are naturally two types of proper
affine spheres, depending on whether the affine normal points toward
or away from the center. For an \emph{elliptic affine sphere}, such
as an ellipsoid, the affine normals point inward toward the center.
\emph{Hyperbolic affine spheres} have affine normals which point
away from the center.  One component of a hyperboloid of two sheets
is the quadric example of a hyperbolic affine sphere.

Some 15 years \c{T}i\c{t}eica's papers, Blaschke's book
\cite{blaschke} records much of the early development of affine
geometry.  Calabi's papers contain many advances in the theory of
affine spheres and related subjects, and Cheng-Yau's resolution of
the structure of hyperbolic affine spheres provides crucial analytic
estimates related to Monge-Amp\`ere equations
\cite{cheng-yau77,cheng-yau86}. We refer the reader to the books of
Nomizu-Sasaki \cite{nomizu-sasaki} and Li-Simon-Zhao
\cite{li-simon-zhao} for overviews of affine differential geometry.
Nomizu-Sasaki \cite{nomizu-sasaki} develop the theory in the modern
notation of connections, while Li-Simon-Zhao \cite{li-simon-zhao}
use Cartan's moving frame techniques and also provide many analytic
details.

This survey article focuses on the relationship between the geometry
of affine spheres to geometric structures on manifolds.  To this end
in Section \ref{2-dim-sec}, we describe the  semi-linear PDE of
\c{T}i\c{t}eica and Wang involving the cubic differential and
developing map for affine spheres in $\re^3$. Then we outline the
relationship, due to Blaschke and Calabi, of affine spheres to real
Monge-Amp\`ere equations and the basic duality results related
related to the Legendre transform in Section \ref{ma-duality}. In
Sections \ref{has-conv-cone}-\ref{proj-man-sec}, we discuss
Cheng-Yau's work on hyperbolic affine spheres and invariants of
convex cones. In Section \ref{affine-manifolds}, we relate affine
spheres to the geometry of affine manifolds, and the conjecture of
Strominger-Yau-Zaslow, relating parabolic affine spheres to
Calabi-Yau manifolds.  In the last two sections, we discuss two
generalizations of affine spheres: Affine maximal hypersurfaces are
generalizations of parabolic affine spheres which have led to
significant progress in the theory of fourth-order elliptic
equations in the solution of Chern's conjecture for affine maximal
surfaces in $\re^3$ by Trudinger-Wang \cite{trudinger-wang00}.  We
also briefly discuss the affine normal flow, which is the natural
parabolic analog of the elliptic PDEs of affine spheres.  The
selection of topics reflects the author's interests, and there are
many important subjects which lie outside the author's expertise. We
largely do not treat nonconvex affine spheres, and we only mention a
few of the recent results classifying affine hypersurfaces with
extremal geometric conditions.

This article is dedicated to Prof.\ S.T.\ Yau, who introduced me to
affine differential geometry, for his rich insight and kind
encouragement over the years, on the occasion of his 59$^{\rm th}$
birthday.

\section{Affine structure equations}

Blaschke used the invariance of the affine normal to derive other
affine-invariant quantities, such as the affine metric and cubic
form \cite{blaschke}.  As in e.g.\ Nomizu-Sasaki
\cite{nomizu-sasaki}, the invariance of the affine normal can be
demonstrated by putting affine-invariant conditions on arbitrary
transverse vector fields to hypersurfaces. So let $L$ be a smooth
strictly convex hypersurface in $\re^{n+1}$ and let $\bar \xi$ be a
transverse vector field. Then we have the following structure
equations of Gauss and Weingarten:
\begin{eqnarray*}
 D_X Y &=& \bar \na_X Y + \bar h(X,Y) \bar \xi, \\
 D_X \bar \xi &=& -\bar S(X) + \bar \tau(X)\bar\xi,
\end{eqnarray*}
where $X,Y$ are tangent vector fields to $L$, $D$ is the standard
affine connection on $\re^{n+1}$, and the equations are given by the
splitting at each $x\in L$ $$ T_x\re^{n+1} = T_x L + \langle \bar
\xi \rangle$$ of the tangent space to $\re^{n+1}$ into the tangent
space to $L$ and the span of $\bar\xi$. In this formulation, $\bar
\na$ is a torsion-free connection on $L$, $\bar h$ is a symmetric
tensor, $\bar S$ is an endomorphism of the tangent bundle $TL$ and
$\bar \tau$ is a one-form.

The affine normal $\xi$  to a convex hypersurface $H$ is the unique
transverse vector field satisfying
\begin{itemize}
\item $\xi$ points to the convex side of $L$ (this is equivalent to
$h$ being positive-definite).
\item $\xi$ is \emph{equiaffine}, which means that $\tau=0$.
\item For $X_i$ a frame of the tangent bundle $TL$, $$ \det_{1\le
i,j \le n} h(X_i,X_j) = \det (X_1,\dots,X_n,\xi)^2,$$
\end{itemize}
where the second determinant is that on $\re^{n+1}$. This approach
is similar to that of \cite{nomizu-sasaki}.

\begin{prop}
The affine normal is well-defined on any smooth strictly convex
hypersurface $H\subset\re^{n+1}$.
\end{prop}
\begin{proof}
Consider any transversal vector field $\bar \xi$ which points to the
convex side of $L$. If we set $$\xi = \phi\bar \xi + Z$$ for $\phi$
a positive function on $L$ and $Z$ a tangent vector field, then
compute, using a frame $X_i$ of the tangent space of $L$, that
 \begin{eqnarray} \phi
&=& \left( \frac{ \det \bar h (X_i,X_j) } {\det (X_1,\dots, X_n,
\bar
\xi)^2 } \right) ^{\frac1{n+2}},\\
 Z^j &=& - \bar h^{ij}(X_i\phi + \phi \bar \tau_i),
\end{eqnarray}
for $\bar h^{ij}$ the inverse of $\bar h_{ij}$.
\end{proof}

\begin{rem}
Given a choice of orientation, the affine normal can be defined on
any $C^3$ nonconvex hypersurface as long as the second fundamental
form is nondegenerate.
\end{rem}

Following Blaschke \cite{blaschke}, we can use the affine normal to
define other affine invariants on $L$.  In particular, we have the
Gauss and Weingarten formulas:
\begin{eqnarray}
\label{gauss-eq}
 D_XY &=& \na_X Y + h(X,Y)\xi, \\
 \label{weingarten-eq}
 D_X\xi &=& -S(X).
\end{eqnarray}
Here $\na$ is the induced \emph{Blaschke connection}, $h$ is the
\emph{affine metric}, or \emph{affine second fundamental form}, and
$S$ is variously called the \emph{affine shape operator}, the
\emph{affine third fundamental form}, or the \emph{affine
curvature}.  The \emph{affine mean curvature} is given by $H =
\frac1n \mbox{tr} S$.

The cubic tensor, given by $C = \na - \hat \na$ for $\hat\na$ the
Levi-Civita connection of $h$, is another important invariant.  Its
main properties are the following:
\begin{prop}
\begin{enumerate}
\item Apolarity: The trace tr\,$ C = 0$ (in index notation, $C^i_{ij}=0$).
\item Symmetry: $C_{ijk} = h_{il}C^l_{jk}$ is totally symmetric on
all three indices.
\item If the cubic form vanishes identically on $L$, $L$ is an open
subset of a hyperquadric.
\end{enumerate}
\end{prop}
The last item is due to Maschke (for analytic surfaces), Pick (for
all surfaces), and Berwald (in general).  See e.g.\
\cite{nomizu-sasaki}.

The symmetry of the cubic form is equivalent to $\na h$ being
totally symmetric.  The Ricci tensor of the affine metric on an
affine sphere is of the form
 \begin{equation} \label{ricci-formula} R_{ij} = (n-1) H g_{ij} +
C_i^{k\ell}C_{jk\ell}. \end{equation} Thus the Ricci tensor is
always bounded from below by $(n-1)H$. This lower bound  is
essential in applying the maximum principle on complete affine
manifolds.

As mentioned above, a proper affine sphere is a hypersurface whose
affine normal lines all converge to a single point, the
\emph{center}.  An improper affine sphere is a hypersurface whose
affine normals are all parallel.  There is an alternate definition
in terms of the affine shape operator:  An \emph{affine sphere} is a
hypersurface whose shape operator is a multiple of the identity $S=
H\, I$, where $H$ is the affine mean curvature. Integrability
conditions then force $H$ to be a constant. The sign of $H$
determines the type of affine sphere: If $H>0$, $L$ is an elliptic
affine sphere. If $H=0$, $L$ is a parabolic affine sphere, and, if
$H<0$, $L$ is a hyperbolic affine sphere.  It is often convenient to
scale proper affine spheres in $\re^{n+1}$ to make the affine mean
curvature $H=\pm1$.

\section{Examples}
The prime examples of affine spheres are the quadric hypersurfaces.
An ellipsoid has affine metric of constant positive curvature.  This
is to be expected, as the ellipsoid is affinely equivalent to a
Euclidean sphere, and the isometries of the Euclidean sphere pull
back to affine actions on the ellipsoid. Similarly, the affine
metric on the hyperboloid has constant negative curvature,
reflecting the Lorentz group action.  An elliptic paraboloid in
$\re^{n+1}$ admits a flat affine metric, and in fact the group of
affine actions preserving the paraboloid is isomorphic to the group
of isometries of $\re^n$.

As we will discuss below, quadrics are the only global examples of
elliptic and parabolic affine spheres.  There are many global
hyperbolic affine spheres, asymptotic to each regular convex cone in
$\re^{n+1}$ (a regular convex cone is an open convex cone which
contains no lines) \cite{calabi72,cheng-yau77,cheng-yau86}.
\c{T}i\c{t}eica already produced the example $$\{(x^1,x^2,x^3) :
x^1x^2x^3=c>0, x^i>0\}$$ in $\re^3$, which is a hyperbolic affine
sphere asymptotic to the boundary of the cone consisting of the
first octant.  Calabi \cite{calabi72} showed the corresponding
example in $\re^{n+1}$ is a hyperbolic affine sphere.
\begin{equation} \label{calabi-example}\left\{(x^1,\dots, x^{n+1}) :
\prod_{i=1}^{n+1} x^i = c, \quad x^i>0 \right\}. \end{equation} The
affine metric of Calabi's example is flat, and its cubic form never
vanishes.

Calabi constructs his example via a product construction for
hyperbolic affine spheres: If $L'\subset \re^{p+1}$ and $L''\subset
\re^{q+1}$ are hyperbolic affine spheres centered at the origin then
the set $$ L = \left\{ \left(x' e^{\frac{-t}{p+1}} , x'' e^{\frac
t{q+1}} \right) : x'\in L', x''\in L'',t\in\re \right\}$$ is a
hyperbolic affine sphere in $\re^{p+q+2}$ (though the affine mean
curvature of $L$ is scaled by a  complicated constant). Applying
this product construction repeatedly to the zero-dimensional
hyperbolic affine sphere $\{1\}\subset\re$ leads to Calabi's example
in the first orthant of $\re^{n+1}$.

We also remark here that in many special cases, geometric conditions
can be imposed on affine spheres (whether convex or not) to
characterize specific equations. There is by now quite a large body
of literature along these lines.  For example, part III of
Nomizu-Sasaki \cite{nomizu-sasaki} details some examples, of which
we mention a few.  Magid-Ryan \cite{magid-ryan} classify all flat
affine spheres (convex or not) in $\re^3$. Convex affine spheres
with whose affine metrics have constant sectional curvature are
shown to be quadrics or affine images of Calabi's example in
Vrancken-Li-Simon \cite{vrancken-li-simon91}. An analogous question
in the non-convex case was settled by Vrancken \cite{vrancken00}.

\section{Two-dimensional affine spheres and \c{T}i\c{t}eica's equation}
\label{2-dim-sec}

\c{T}i\c{t}eica first studied affine spheres in $\re^3$ (more
properly, he studied a subset of affine spheres he called
$S$-surfaces) \cite{tzitzeica08,tzitzeica09}, and found conditions
under which these surfaces can be integrated from initial data. More
specifically, if $\alpha,\beta$ are two real parameters, and $v =
v(\alpha,\beta)$ satisfies \begin{equation} \label{tzitz-eq}
 \frac{\partial ^2 v} {\partial \alpha \partial \beta} = e^v -
 e^{-2v},
\end{equation}
then the system of equations
\begin{eqnarray} \label{tz1}
\frac{\partial ^2 f}{\partial \alpha^2} &=& \frac{\partial
v}{\partial \alpha} \frac{\partial f}{\partial \alpha} + e^{-v}
\frac{\partial f}{\partial \beta}, \\
\label{tz2} \frac{\partial ^2 f}{\partial \beta^2} &=& e^{-v}
\frac{\partial f}{\partial \alpha} + \frac{\partial v}{\partial
\beta}
\frac{\partial f}{\partial \beta}, \\
\label{tz3} \frac{\partial ^2f}{\partial \alpha \partial \beta} &=&
e^v f
\end{eqnarray}
is integrable for $f=f(\alpha,\beta)$ a map into $\re^3$. This
system can considered as a first-order system in the frame
$\{f,\frac{\partial f}{\partial \alpha}, \frac{\partial f}{\partial
\beta}\}$, and \c{T}i\c{t}eica's equation (\ref{tzitz-eq}) is the
integrability condition.  More specifically, if $f\!:\mathcal
D\to\re^3$ from a simply connected domain $\mathcal D\subset\re^2$,
we may specify $f(x_0), f_\alpha(x_0), f_{\beta}(x_0)$ for any point
$x_0\in \mathcal D$.  Equations (\ref{tz1}-\ref{tz3}) can be
considered as a first-order system of PDEs in
$\{f,f_\alpha,f_\beta\}$, and thus can be integrated along any path
from $x_0$.  Then (\ref{tzitz-eq}) shows that the solution at any
$x\in \mathcal D$ is independent of the path chosen from $x_0$ to
$x$.  The integrability condition is determined by checking  e.g.\
$(f_{\alpha\beta})_ \alpha = (f_{\alpha\alpha})_\beta$ and using the
Frobenius Theorem.

In modern language, the surface parametrized by $f$ is then an
improper, nonconvex affine sphere centered at the origin. This means
the affine metric is indefinite. But as is often the case in the
theory of two-dimensional integrable systems, the signature of the
metric can be changed by considering complex parameters. So we can
produce convex as well as nonconvex affine spheres by
\c{T}i\c{t}eica's method.  In fact, \c{T}i\c{t}eica produces the
affine sphere
$$\{(x^1,x^2,x^3) : x^1 x^2 x^3 = 1, x^i>0 \}.$$

Wang \cite{wang91} and Simon-Wang \cite{simon-wang93} have extended
\c{T}i\c{t}eica's technique to all convex affine spheres in $\re^3$.
The affine metric is positive definite on a convex surface in
$\re^3$, and so there is an induced conformal structure.  So in this
case, we may assume that $f\!: \mathcal D\to \re^3$, where $\mathcal
D\subset \co$ is simply connected and so that the map is a conformal
map with respect to the affine metric on the image $f(\mathcal D)$.
Choose a local complex coordinate $z$ on $\mathcal D$, so that the
affine metric is $$e^\psi|dz|^2.$$ Then the apolarity condition on
the cubic form shows that all but two components of the cubic form
vanish.  In terms of complex coordinates, if we set $$ U = C^{\bar
1} _{11} e^\psi,$$ then the structure equations
(\ref{gauss-eq}-\ref{weingarten-eq}) for the (complexified) frame
$\{f_z, f_{\bar z} ,\xi\}$ become
\begin{eqnarray*}
f_{zz}&=& \psi_z f_z + U e^{-\psi} f_{\bar z}, \\
f_{\bar z \bar z} &=& \bar U e^{-\psi} f_z + \psi_{\bar z} f_{\bar
z}, \\
f_{z \bar z} &=& \sfrac12 e^\psi \xi,\\
\xi_z &=& -H f_z, \\
\xi_{\bar z} &=& -H f_{\bar z}.
\end{eqnarray*}
The integrability conditions for these equations are then
\begin{eqnarray}
\label{U-holomorphic}
 U_{\bar z} &=& 0, \\
 \label{local-tz-eq}
 \psi_{z \bar z} + |U|^2 e^{-2\psi} + \frac H 2 e^\psi &=& 0.
\end{eqnarray}

If the holomorphic coordinate $z$ is changed, $U$ transforms as a
cubic differential, which is holomorphic by (\ref{U-holomorphic}).
On a Riemann surface, (\ref{U-holomorphic}) becomes for $e^\psi
|dz|^2 = e^\phi g$, \begin{equation} \label{as-eq-2} \Delta \phi +
4\|U\|^2 e^{-2\phi} + 2H e^\phi - 2\kappa,
 \end{equation}
where $\Delta$ is the Laplace operator of $g$, $\|\cdot\|$ is the
induced norm on cubic differential, and $\kappa$ is the Gauss
curvature. Below we will discuss global solutions to (\ref{as-eq-2})
on Riemann surfaces due to
\cite{labourie97,labourie07,loftin02c,loftin04, lyz05,lyz-erratum,
wang91}, and their application to the geometry of projective,
affine, and Calabi-Yau manifolds.

(We note different versions of \c{T}i\c{t}eica's equation also occur
in other geometric contexts in which real forms of
$\mathbf{SL}(3,\co)$ act. For example, see McIntosh
\cite{mcintosh03} for an application  of solutions to $$ \psi_{z
\bar z} - e^{-2\psi} + e^\psi = 0$$ to minimal Lagrangian immersions
in $\cp^2$ and special Lagrangian cones in $\co^3$.)

The structure equations (\ref{gauss-eq}-\ref{weingarten-eq}) are a
first-order linear PDE system in the frame $\{\xi,f_z,f_{\bar z}
\}$, and if (\ref{U-holomorphic} - \ref{local-tz-eq}) are satisfied,
then (\ref{gauss-eq}-\ref{weingarten-eq}) can be solved as an
initial value problem on any simply connected domain $\mathcal
D\subset \co$.  In other words, for any $z_0\in \mathcal D$, if
$\xi(z_0),f_z(z_0),f_{\bar z}(z_0)$ are specified, then
(\ref{gauss-eq}-\ref{weingarten-eq}) determine the frame at every
point in $\mathcal D$.

For affine spheres, this determines an affine or projective holonomy
action.  Proper affine spheres with center at the origin naturally
have holonomy in $\Sl{3}$ (affine actions in $\sa{3}$ fixing the
origin), while improper affine spheres with affine normal $(0,1)$
naturally have holonomy in $\sa 2$ (actions in $\sa{3}$ fixing the
affine normal).  For the special case of a holomorphic equivalence
relation $z\sim z+1$, the frame $\{\xi, f_z ,f_{\bar z}\}$ is a
well-defined frame on the Riemann surface $\mathcal D /\sim$.  For
example, for a hyperbolic affine sphere with center 0 and affine
mean curvature $-1$, the affine normal $\xi = f$, and integrating
the structure equations along a path from $z_0$ to $z_0 + 1$
calculates holonomy map in $\Sl3$ for the frame $\{f,f_z,f_{\bar z}
\}$. In particular, on a Riemann surface with a point singularity,
it is often possible to prescribe the behavior of the affine metric
and cubic differential near the singularity, and then to use the
theory of ODEs to determine the conjugacy class of holonomy of a
loop around the singularity. See \cite{loftin02c,loftin04}.

Parabolic affine spheres in $\re^3$ have much better integrability
properties.  Much like minimal surfaces in $\re^3$, there are
Weierstrass formulas for reproducing parabolic affine spheres in
terms of holomorphic functions.

The basic ideas leading to this Weierstrass formula have been
available for quite a long time.  Blaschke \cite{blaschke}
recognizes that parabolic affine spheres with affine normal
$(0,0,1)$ are locally given in terms of a graph of a convex function
$u$ satisfying the Monge-Amp\`ere equation \begin{equation}
\label{ma-2dim} \det \frac{\partial ^2 u }{ \partial x^i \partial
x^j} = 1,
\end{equation}
and he shows that parabolic affine spheres in $\re^3$ are completely
integrable, by using the observation in Darboux \cite{darboux} that
solutions to (\ref{ma-2dim}) in dimension two can locally be
transformed into harmonic functions. J\"orgens \cite{jorgens54} uses
the relation between solutions to (\ref{ma-2dim}) and complex
analytic functions to show that any entire convex solution to
(\ref{ma-2dim}) in $\re^2$ must be a quadratic polynomial. Moreover,
J\"orgens uses a natural transformation between parabolic affine
two-spheres and minimal surfaces to reprove Bernstein's Theorem that
any minimal surface in $\re^3$ which is an entire graph is a plane.
See also Chapter 9 of Spivak \cite{spivak-vol4}.

The affine Weierstrass formula for affine maximal surfaces in
$\re^3$ (a generalization of a parabolic affine sphere) is given in
by Calabi \cite{calabi88}, Terng \cite{terng83}, and Li \cite{li89},
for affine maximal surfaces in $\re^3$. A similar description for
parabolic affine spheres is given later by Ferrer-Mart\'inez-Mil\'an
\cite{ferrer-mm99}, although an equivalent formulation motivated by
string theory is found earlier in Greene-Shapere-Vafa-Yau
\cite{greene-svy90}. The affine Weierstrass formula of
\cite{ferrer-mm99} is the following: Given two holomorphic functions
$F$ and $G$ satisfying $|dF| < |dG|$ on a simply connected domain
$\mathcal D\subset \co$, the parametrized surface $$ \left( \frac12
(G+\bar F), \frac18 (|G|^2-|F|^2), \frac14 \mbox{Re}\, (GF) -
\frac12 \mbox{Re}\int F\,dG \right)$$ is a parabolic affine sphere
with affine normal $(0,0,1)$.

There is a higher-dimensional generalization of the Weierstrass
representation formula to special K\"ahler manifolds. A special
K\"ahler manifold is a K\"ahler manifold which admits a
torsion-free, flat connection $\na$ with respect to which the
K\"ahler form is parallel and so that $d_{\na}I=0$ for $I$ the
complex structure tensor \cite{freed99}. Cort\'es \cite{cortes02}
has shown that any special K\"ahler manifold may be constructed from
the graph of a holomorphic function on a domain in $\co^n$. Z.\ Lu
has shown that any complete special K\"ahler manifold is flat
\cite{lu99}.  Lu's result also follows from Calabi's Theorem
\ref{calabi-affine-complete} below, since each special K\"ahler
manifold carries the structure of a parabolic affine sphere
\cite{baues-cortes01}.

It is also possible to relate solutions to other versions of
\c{T}i\c{t}eica's equation to integrable systems, as
Dunajski-Plansangkate have recently related radially symmetric
solutions of (\ref{local-tz-eq}) with $H=1$ to solutions of
Painlev\'e III.

\section{Monge-Amp\`ere Equations and Duality} \label{ma-duality}

In this section, we recount the real Monge-Amp\`ere equations
related to convex affine spheres, and use the conormal map and the
Legendre transform to find dual affine spheres.  In this context,
the equation for a parabolic affine sphere, $\det u_{ij} =1$, goes
back to Blaschke. The equations for proper affine spheres are due to
Calabi \cite{calabi72}; although we will primarily present them in
the context of Gigena \cite{gigena81}.

For simplicity, we state the following theorem only for affine
spheres of affine mean curvature $H=-1,0,1$.  More general proper
affine spheres may be obtained by scaling.

\begin{prop}
\begin{itemize}
\item A hyperbolic affine sphere with center 0 and affine mean
curvature $-1$ is locally given by the radial graph of $-1/u$,
$$\left\{ -\frac 1{u(t)} (1,t^1,\dots,t^n) :
t=(t^1,\dots,t^n)\in\Omega.\right\},$$ where $\Omega$ is a domain in
$\re^n$ (thought of as an inhomogeneous domain in $\rp^n$), and $u$
is a convex negative function satisfying
$$\det \frac{\partial ^2u}{\partial t^i \partial t^j} = \left(-\frac1u\right)^{n+2}.$$
The affine metric is given by $$-\frac1u \, \frac{\partial^2
u}{\partial t^i \partial t^j} \, dt^i dt^j.$$
\item A parabolic affine sphere with affine normal $(0,\dots,0,1)$
is given by the graph of a convex function $u$
$$\{(x^1,\dots,x^n,u(x)) : x=(x^1,\dots,x^n)\in \mathcal O \subset
\re^n$$ satisfying the Monge-Amp\`ere equation $$\det u_{ij} = 1.$$
The affine metric is given by $$\frac{\partial ^2u}{\partial x^i
\partial x^j} \, dx^i dx^j.$$
\item An elliptic affine sphere with center 0 and affine mean
curvature $1$ is given by the radial graph of $1/u$, $$\left\{ \frac
1{u(t)} (1,t^1,\dots,t^n) : t=(t^1,\dots,t^n)\in\Omega.\right\},$$
where $\Omega$ is a domain in $\re^n$ (thought of as an
inhomogeneous domain in $\rp^n$), and $u$ is a convex positive
function satisfying
$$\det \frac{\partial ^2u}{\partial t^i \partial t^j} = \left(\frac1u\right)^{n+2}.$$
The affine metric is given by $$ \frac1u \, \frac{\partial ^2
u}{\partial t^i \partial t^j} \, dt^i dt^j.$$
\end{itemize}
\end{prop}

For each convex affine sphere as above in $\re^{n+1}$, there is a
dual affine sphere in $\re_{n+1}$, the dual vector space to
$\re^{n+1}$.   The construction is slightly different in the case of
proper and improper affine spheres, but both constructions are
related to the Legendre transform.  Given a smooth strictly convex
function $v \!:\Omega\to \re$ on a convex domain
$\Omega\subset\re^n$, the Legendre transform function $v^*$ is
defined by $$ v^* + v = x^i \, \frac{\partial v} {\partial x^i}.$$
The function $v^*$ is considered primarily as a convex function $v^*
: \Omega^{*,v} \to \re$ with
 $$\Omega^{*,v} = \left\{ \left(
\frac{\partial v}{\partial x^i}(x) \right) \in \re_n : x\in
\Omega\right\}.
$$

The duality of a parabolic affine sphere comes directly from the
Legendre transform, while for proper affine spheres, the duality is
provided by the conormal map. Given a hypersurface $L\subset
\re^{n+1}$ which is transverse to the position vector, the conormal
map  $N\!:L\to \re_{n+1}$ is given for $x\in L$ by
$$N\!:x\mapsto \ell ,\qquad \ell(x) =1, \quad \ell(T_xL)=0.$$
The conormal map is naturally related to the Legendre transform by
the following formulas: If $u = u(t^1,\dots,t^n)$ is convex, and $L$
is the radial graph of $c/u$: $$L = \{(c/u)(t^1,\dots,t^n,1)\},$$
then the conormal map of $L$ is given by
 \begin{equation} \label{conormal-legendre} \left\{ -\frac1c \left(
-\frac{\partial u}{\partial t^1} , \dots, -\frac{\partial
u}{\partial t^n} , u^* \right) \right\} \end{equation}
 for $u^*$ the
Legendre transform of $u$.

\begin{prop} \label{dual-affine-sphere}
To each affine sphere in $\re^{n+1}$, there is a dual affine sphere
of the same type in the dual space $\re_{n+1}$.  The dual of a
proper affine sphere centered at the origin is given by the image of
the conormal map.

For an improper affine sphere with affine normal $\xi=(0,\dots,0,1)$
given by the graph $(x,u(x))$, the dual is the graph of the Legendre
transform of $u$.
\end{prop}

\begin{rem}
Proposition \ref{dual-affine-sphere} and the relationship between
the conormal map and the Legendre transform
(\ref{conormal-legendre}) lead to Calabi's original formulation of
the relationship between proper affine spheres and solutions to the
Monge-Amp\`ere equation \cite{calabi72}: The graph of a convex
function $$\{(x,\psi(x)): x\in\Omega\subset\re^n\}$$ is a proper
affine sphere of affine mean curvature $H$ and center 0 if and only
if the Legendre transform $\psi^*$ of $\psi$ satisfies $$\det
\psi^*_{ij} = (H\psi^*)^{-n-2}.$$
\end{rem}

We need the notion of conjugate connections to explain the duality
further.  Given a smooth vector bundle $E\to M$, equipped with a
nondegenerate metric $h$ and a connection $\na$, the \emph{conjugate
connection} $\na^*$ on $E$ is given by
$$ d\, h(s,t) = h(\na s, t) + h (s, \na^*t)$$ for smooth sections
$s,t$ of $E$.  In the case $\na h$ is totally symmetric, which is
always true for affine structures, the Levi-Civita connection $\hat
\na$ of $h$ satisfies \begin{equation} \label{conj-levi-civita} \hat
\na = \sfrac12(\na + \na ^*). \end{equation} See e.g.\
\cite{nomizu-sasaki}.

\begin{prop} \label{dual-conj}
The Blaschke connection on an improper affine sphere is flat.  The
Blaschke connection on a proper affine sphere is projectively flat.

In all cases, the dual map is an isometry with respect to the affine
metrics. It maps the cubic form to its negative; in other words, the
Blaschke connections of the dual hyperspheres are conjugate with
respect to the affine metric.
\end{prop}

We sketch the proofs of Propositions \ref{dual-affine-sphere} and
\ref{dual-conj} together.

\begin{proof}[Sketch of proof]
For improper affine spheres, it is well known that the Legendre
transform of a convex solution to the Monge-Amp\`ere equation $\det
u_{ij} = 1$ solves the same equation.  Moreover, the affine metric
for a parabolic affine sphere is of Hessian type $u_{ij} dx^i dx^j$,
and the affine connection is flat. For such metrics, the conjugate
connection is also flat \cite{amari-nagaoka}, and corresponds to the
dual affine structure.

For any convex affine hypersurface $L$ in $\re^{n+1}$, the image of
the conormal map is naturally a centroaffine hypersurface
$L^*\subset \re_{n+1}$ (in other words, $L^*$ is transverse to the
position vector). The centroaffine connection of $L^*$ is conjugate
to the affine (Blaschke) connection of $L$
\cite{schirokov,nomizu-sasaki}.  For $f$ the position vector of
$L^*$, the centroaffine connection $\na^c$ is defined by the Gauss
equation $$ D_X Y = \na^c_X Y + h^c(X,Y)f$$ defined by the splitting
$T_f\re_{n+1} = T_f L^* + \langle f \rangle.$  But since the affine
normal for proper affine spheres with center 0 and affine mean
curvature $\pm1$ is equal to $\mp f$, $L^*$ is again a proper affine
sphere with Blaschke connection $\na^c$.

Note that (\ref{conj-levi-civita}), by the definition of the cubic
form, shows that passing from the Blaschke connection to its
conjugate corresponds to mapping the cubic form to its negative.
\end{proof}

\section{Global classification of affine spheres}

The affine metric of a convex affine hypersurface in $\re^{n+1}$ is
a Riemannian metric, and thus gives an affine-invariant notion of
completeness on the hypersurface.  We may also consider whether the
convex hypersurface is properly embedded, which corresponds to an
extrinsic notion of \emph{Euclidean completeness} for hypersurfaces.
In general, the two notions of completeness are different (see e.g.\
examples in \cite{li-simon-zhao}), but for affine spheres, the two
notions are the same. This was conjectured by Calabi and proved, in
the hyperbolic case, by Cheng-Yau.

\begin{thm}[J\"orgens, Calabi, Pogorelov, Cheng-Yau]
Let $\Omega$ be a domain in $\re^n$ and let $u\!: \Omega\to \re$ be
a smooth convex function satisfying the Monge-Amp\`ere equation
$$\det u_{ij} = 1$$ and so that the graph
$$\{(x,u(x)) : x\in \Omega\}$$ is closed in $\re^{n+1}$.  Then $u$
is a quadratic polynomial.
\end{thm}

This theorem is due to J\"orgens originally for entire solutions
when $n=2$ \cite{jorgens54}, to Calabi for entire solutions for
$n=3,4,5$ \cite{calabi58}, and Pogorelov for entire solutions for
general $n$ \cite{pogorelov72}.  Cheng-Yau proved the full theorem
in \cite{cheng-yau86}. In terms of affine geometry,  this theorem is
equivalent to the following:
\begin{cor}
Every convex properly embedded parabolic affine sphere in
$\re^{n+1}$ is an elliptic paraboloid.
\end{cor}

The corresponding theorem for affine-complete parabolic affine
spheres is due to J\"orgens in dimension two \cite{jorgens54} and
Calabi in general \cite{calabi58}:

\begin{thm} \label{calabi-affine-complete}
Any affine-complete parabolic affine sphere in $\re^{n+1}$ is an
elliptic paraboloid.
\end{thm}

The proof is to use a maximum principle on noncompact manifolds to
show that the norm of the cubic form must vanish. Then Berwald's
theorem shows that the affine sphere must be a quadric hypersurface.

Similarly, all global examples of elliptic affine spheres must be
ellipsoids. Blaschke first shows that any compact elliptic affine
sphere in $\re^3$ is an ellipsoid \cite{blaschke}, and Deicke later
extends this theorem to any dimension \cite{deicke53}.

Any affine-complete elliptic affine sphere must also be an
ellipsoid.  Calabi \cite{calabi72} shows that the Ricci curvature of
an elliptic affine sphere is positive, and thus Myers's Theorem
shows the affine sphere must be compact.  Euclidean-complete
elliptic affine spheres are affine complete by estimates of
Cheng-Yau \cite{cheng-yau86}.  Thus any affine or Euclidean complete
affine sphere in $\re^{n+1}$ must be compact, and thus is an
ellipsoid.

For any strictly convex smooth affine hypersurface in $\re^{n+1}$,
Trudinger-Wang prove that affine completeness implies Euclidean
completeness if $n\ge 2$
\cite{trudinger-wang02}.


\section{Hyperbolic affine spheres and invariants of convex cones}
\label{has-conv-cone}

Let $\re_{n+1}$ denote the dual vector space of $\re^{n+1}$. Then
given an open convex cone $\mathcal C\subset \re^{n+1}$, the dual
cone $\mathcal C^* \subset \re_{n+1}$ can be defined as $$\mathcal
C^* = \{ \ell \in \re_{n+1} : \ell(x) > 0 \mbox{ for all } x\in
\mathcal C \}.$$

The existence of hyperbolic affine spheres is due to Cheng-Yau. In
\cite{cheng-yau77}, they show that for any convex bounded domain
$\Omega\subset\re^n$, there is a unique convex solution to the
Dirichlet problem
 \begin{equation} \label{has-ma} \det u_{ij} = \left(-\frac1u\right)^{n+2}, \qquad
u= 0 \mbox{ on }\partial \Omega. \end{equation} Let $$\mathcal C  =
\{ s(1,t) : t\in\Omega, s>0\}$$ be the cone over $\Omega$.  Then
they show in \cite{cheng-yau86} that $u$ induces the hyperbolic
affine sphere asymptotic to the boundary of the dual cone $\mathcal
C^*$ by taking the Euclidean graph of the Legendre transform of $u$.
The description of the hyperbolic affine sphere asymptotic to
$\mathcal C$ is given in Gigena \cite{gigena81} as the radial graph
of $-1/u$.

We outline a proof of Cheng-Yau's solution to the Monge-Amp\`ere
equation (\ref{has-ma}).  By now, producing convex solutions to
Dirichlet boundary problems for real Monge-Amp\`ere equations $\det
v_{ij} = F(x,v)$ on strictly convex bounded domains for smooth
positive functions $F$ is fairly standard (see e.g.\ \cite[Theorem
17.22]{gilbarg-trudinger}). Equation (\ref{has-ma}) is singular in
two ways, however.  First of all, since $u=0$ on $\partial \Omega$,
the right-hand side $(-\frac1u)^{n+2}$ blows up on $\partial
\Omega$.  Second, $\Omega$ is allowed to be any convex domain with
potentially no more than Lipschitz boundary regularity. Below we
give limiting arguments, essentially due to  Cheng-Yau
\cite{cheng-yau77}, to produce solutions to (\ref{has-ma}) in the
general case. Calabi's example (\ref{calabi-example}) above provides
an explicit solution to (\ref{has-ma}) on a  simplex, which is used
as a barrier.  In dimension two, Loewner-Nirenberg
\cite{loewner-nirenberg} solved (\ref{has-ma}) on domains with
strictly convex smooth boundary, and noted the projective invariance
of solutions to (\ref{has-ma}) from a point of view independent of
affine geometry.

\begin{proof}[Sketch of proof]
Let $\Omega\subset\re^n$ be a bounded convex domain, and consider an
exhaustion $$\Omega = \bigcup_k \Omega_k, \qquad \Omega_k \Subset
\Omega_{k+1},$$ for $\Omega_k$ strictly convex bounded domains with
smooth boundary.  Standard theory for the real Monge-Amp\`ere
equation shows that there is a unique convex solution
$u_{k,\epsilon}$ to \begin{equation} \label{perturbed-has-eq} \det
u_{ij} = \left(-\frac1 {u-\epsilon} \right)^{n+2} \mbox{ on }
\Omega_k, \qquad u|_{\partial \Omega_k} = 0
\end{equation}
 for each $\epsilon>0$. Each $u_{k,\epsilon}$ is smooth on
$\Omega_k$ and is continuous on $\bar \Omega_k$.  Our solution $u$
will be the limit of $u_{k,\epsilon}$ as $k\to\infty$ and
$\epsilon\to0$.

First, we let $\epsilon\to 0$. The maximum principle shows that
$u_{k,\epsilon} \le u_{k,\epsilon'}$ if $\epsilon < \epsilon'$: We
will find a contradiction if the difference $w=u_{k,\epsilon} -
u_{k,\epsilon'}$ is positive anywhere on $\bar\Omega_k$.  Since $w$
vanishes on $\partial\Omega_k$, if $w$ is positive anywhere, it has
a positive maximum point $p\in\Omega_k$.  Then at $p$, the Hessian
$$ \frac{\partial ^2 w}{\partial x^i \partial x^j}(p) = \frac{\partial ^2
u_{k,\epsilon}}{\partial x^i \partial x^j}(p) -\frac{\partial ^2
u_{k,\epsilon'}}{\partial x^i \partial x^j}(p)$$ is negative
definite.  Since the Hessians of both $u_{k,\epsilon}$ and
$u_{k,\epsilon'}$ are positive-definite, a simple lemma shows that
$$\det \left( \frac{\partial ^2
u_{k,\epsilon}}{\partial x^i \partial x^j}(p) \right) \le \det
\left( \frac{\partial ^2 u_{k,\epsilon'}}{\partial x^i \partial
x^j}(p) \right) .$$ But then, the equation
(\ref{perturbed-has-eq}) shows $$\left(-\frac1
{u_{k,\epsilon}-\epsilon}(p) \right)^{n+2} \le \left(-\frac1
{u_{k,\epsilon'}(p)-\epsilon'} \right)^{n+2}.$$ This
contradicts the assumptions that $\epsilon<\epsilon'$ and $w(p)
= u_{k,\epsilon}(p) - u_{k,\epsilon'}(p)>0$.

Therefore, as $\epsilon\to0$, $u_{k,\epsilon}$ is decreasing
pointwise, and thus there is a pointwise limit function
$u_k=u_{k,0}$ unless the sequence decreases to $-\infty$. By
using affine transformations, Calabi's example provides an
explicit solution $v$ to (\ref{has-ma}) on any simplex
$\Delta\subset\re^n$. For a given boundary point $q\in\partial
\Omega_k$, choose $\Delta$ so that $\Delta\supset \Omega_k$ and
$q\in \partial \Delta \cap \partial \Omega_k$. The maximum
principle then shows that $u_{k,\epsilon}\ge v$ on all of $\bar
\Omega_k$, and so the limit function does not go to $-\infty$
and is continuous at the boundary point $q$. Estimates for the
Monge-Amp\`ere equation show the convergence $u_{k,\epsilon}
\to u_k$ is in $C^\infty_{\rm loc} (\Omega_k)$ and thus $u_k$
solves (\ref{has-ma}): $C^1$ estimates are standard for convex
functions, while Pogorelov provides interior $C^2$ estimates
\cite{pogorelov72}, and then one can use either the interior
$C^3$ estimates of Calabi \cite{calabi58} or interior
$C^{2,\alpha}$ estimates of Evans to achieve the desired
regularity.

Now $u_k$ solves (\ref{has-ma}) on $\Omega_k$. Now the same
ideas allow us to take $k\to\infty$: The maximum principle
implies $u_{k+1} \le u_k$ on $\Omega_k$ and Calabi's barriers
ensures the limit $u_k\to u$ is finite and continuous to the
boundary of $\Omega$.  Interior estimates as mentioned above
implies the convergence $u_k\to u$ is in $C^\infty_{\rm
loc}(\Omega)$, and so $u$ solves (\ref{has-ma}) on $\Omega$.
\end{proof}

The following theorem is due to Cheng-Yau
\cite{cheng-yau77,cheng-yau86} and Calabi-Nirenberg (unpublished),
with clarifications due to Gigena \cite{gigena81}, Sasaki
\cite{sasaki80} and A.-M.\ Li \cite{li90,li92}.

\begin{thm}
For any open convex cone $\mathcal C \subset \re^{n+1}$ which
contains no lines, there is a unique convex properly embedded
hyperbolic affine sphere $L\subset\re^{n+1}$ which has affine mean
curvature $-1$, has the vertex of $\mathcal C$ as its center, and is
asymptotic to the boundary $\partial C$.  For any immersed
hyperbolic affine sphere $L\to\re^{n+1}$, the properness of the
immersion is equivalent to the completeness of the affine metric,
and any such $L$ is a properly embedded hypersurface asymptotic to
the boundary of the cone $\mathcal C$ given by the convex hull of
$H$ and its center.

The image of the conormal map of $L$ is the unique hyperbolic affine
sphere of affine mean curvature $-1$ vertex 0 asymptotic in the dual
space $\re_{n+1}$ to the boundary of the dual cone $\mathcal C^*$.
\end{thm}

\begin{proof}[Sketch of proof]
Cheng and Yau \cite{cheng-yau86} prove (with a small gap) that an
affine sphere is properly embedded if and only if its affine metric
is complete.  Moreover, any such hyperbolic affine sphere $L$ is
asymptotic to the cone given the convex hull of $L$ and its center.
(Calabi-Nirenberg, in unpublished work, establish the same result.)
In \cite{li90,li92}, A.-M.\ Li clarified the proof of Cheng-Yau by
using essentially the same estimates developed in \cite{cheng-yau86}
to show that affine completeness implies Euclidean completeness for
hyperbolic affine spheres.  Trudinger-Wang have recently proved that
any convex affine-complete hypersurface in $\re^{n+1}$ is Euclidean
complete if $n\ge2$ \cite{trudinger-wang02}.

Given a bounded convex domain $\Omega\subset\re^n$, considered as
lying in an affine subspace $\{x^0 = 1 \} \subset \re^{n+1}$, denote
the cone over $\Omega$ as $$\mathcal C =  \bigcup_{s>0}
\left\{s(1,t) : t\in \Omega \right\}.$$  Consider Cheng-Yau's
solution to the Monge-Amp\`ere equation $\det
u_{ij}=(-\frac1u)^{n+2}$ with zero Dirichlet boundary values
\cite{cheng-yau77}. The hyperbolic affine sphere $L$ in the cone
$\mathcal C\subset \re^{n+1}$ given by the radial graph of $-1/u$ is
then obviously asymptotic to the boundary of $\mathcal C$ and is
Euclidean complete, by work of Gigena \cite{gigena81}. (This
formulation was known experts in the 1970s, when Cheng-Yau completed
their work.) An alternate approach is taken by Sasaki
\cite{sasaki80}, who shows that the dual affine sphere to $L$, given
by the Legendre transform of $u$, is a Euclidean complete hyperbolic
affine sphere $L^*$ asymptotic to the dual cone $\mathcal C^*
\subset \re_{n+1}$.

These two affine spheres $L$ and $L^*$ are dual under the conormal
map. The local statement of duality follows from the identification
of the centroaffine connection of the image of the conormal map with
the conjugate to the Blaschke connection, already found in the
Schirokovs' book \cite{schirokov}. In other words, the dual of $L$
is a hyperbolic affine sphere.  To show that the image of the
conormal map of $L$ is $L^*$ follows from  Gigena \cite{gigena81}.
Alternately, once we produce either $L$ or $L^*$ as a
Euclidean-complete affine sphere by \cite{gigena81} or
\cite{sasaki80}, it must be affine complete. Since the conormal map
is an isometry of the affine metric, its dual $L^*$ or $L$
respectively is also affine complete, and so it must be Euclidean
complete.  We can identify this dual with the correct target by the
uniqueness of solutions to (\ref{has-ma}), which follows from the
maximum principle.

We also indicate very briefly some of Cheng-Yau's estimates
\cite{cheng-yau86}:  To show a Euclidean-complete affine sphere $L$
is affine-complete, first of all note that we can choose affine
coordinates on $\re^{n+1}$ so that one coordinate function (a height
function $\mathcal H$) is proper on $L$.  Then on the compact domain
$L_c = \{\mathcal H \le c\}$, use the maximum principle to estimate
functions of the form $$ \exp \left( \frac{-m}{c-\mathcal H} \right)
\frac {|\hat \na \mathcal H|^2}{(\mathcal H + \epsilon )^p}
$$ for appropriate constants $m,\epsilon,p$ and $\hat\na$ the Levi-Civita
connection of the affine metric.  Then, taking $c\to\infty$ produces
an estimate of the form $$\frac{|\hat\na \mathcal H  |} { \mathcal H
+1 } \le C$$ for a constant $C$.  Since $\log(\mathcal H + 1)$ is
proper on $L$, this gradient estimate shows the affine metric is
complete.

To prove that affine completeness implies Euclidean completeness for
hyperbolic affine spheres, a gradient estimate for the Legendre
transform of $\mathcal H$ \cite{cheng-yau86,li90,li92} and an
estimate of the norm of the cubic form \cite{calabi72} are used.
\end{proof}

Calabi shows that the affine metric on complete hyperbolic affine
spheres has Ricci curvature bounded between 0 and a negative
constant \cite{calabi72}.  The lower bound is a  pointwise formula
for the Ricci tensor (\ref{ricci-formula}), while the fact that the
Ricci is nonpositive requires global techniques from Riemannian
geometry, and a bound on the norm of the cubic form.  We note that
the extreme cases are both satisfied by homogeneous affine spheres:
The affine metric on a hyperboloid has constant negative sectional
curvature, while the affine metric on Calabi's example
(\ref{calabi-example}) is flat.  It is instructive to think of these
examples in terms of extrema of convex projective domains. The
hyperboloid is asymptotic to a cone over a round ball, while
Calabi's example is asymptotic to a cone over a simplex.

The Monge-Amp\`ere equation (\ref{has-ma}) is  very similar to
Fefferman's equation for complete K\"ahler-Einstein metrics on
bounded smooth strictly pseudoconvex domains in $\co^n$
\cite{sasaki85}. Sasaki has used a similar asymptotic calculations
to Fefferman's to compute projective invariants of convex domains
\cite{sasaki88}.

\section{Projective Manifolds} \label{proj-man-sec}
An $\rp^n$ structure (real projective structure) on a smooth
manifold $M$ is given by an atlas of coordinate charts in $\rp^n$
with gluing maps locally constant projective maps in $\pgl{n+1}$.
(So in Thurston's language, an $\rp^n$ manifold is an $(X,G)$
manifold for the homogeneous space $X=\rp^n$ with group
$G=\pgl{n+1}$.)  An $\rp^n$ manifold $M$ is \emph{properly convex}
if it is given by a quotient
$$M = \Omega / \Gamma,$$ where $\Omega$ is a bounded convex domain
in $\re^n \subset \rp^n$ and $\Gamma\subset \pgl{n+1}$ acts
discretely and properly discontinuously on $\Omega$.  Hyperbolic
manifolds all admit $\rp^n$ structures by the Klein model of
hyperbolic space: In this model hyperbolic space $\mathbb H^n$ is
represented by an open ball $B \subset \re^n$, and the hyperbolic
isometries are exactly the projective actions on $\rp^n\supset
\re^n$ which act on $B$.

A \emph{geodesic} in an $\rp^n$ manifold is a path which is a
straight line segment in each $\rp^n$ coordinate chart.  This leads
to an alternate definition of an $\rp^n$ structure: An $\rp^n$
structure on $M$ is given by a projective equivalence class of
projectively flat connections on the tangent bundle $TM$.  Two
connections on $TM$ are said to be \emph{projectively equivalent} if
they have the same geodesics, when considered as unparametrized
sets.  A connection is \emph{projectively flat} if it is locally
projectively equivalent to a flat, torsion-free connection.  The
relation between the two definitions of $\rp^n$ manifold is the
following: The geodesics of the $\rp^n$ structure are then exactly
the geodesics of each projectively flat connection in the
equivalence class.  (Yet another equivalent definition of an $\rp^n$
structure on a manifold $M$ is the existence of a flat projective
\emph{Cartan connection} on $M$. See e.g.\ Kobayashi
\cite{kobayashi72}.)

For a hyperbolic affine sphere $L$ with center 0 and affine mean
curvature $-1$, Gauss's equation reads
$$ D_XY = \na_XY + h(X,Y) f,$$
where $f$ is the position vector.  In particular, in this case, the
geodesics of $\na$ project to straight lines under the projection
$L\to\Omega$ (see e.g.\ \cite{nomizu-sasaki}).  This shows that
$\na$ is projectively flat.  Note this argument shows that any
centroaffine connection $\tilde \na$ is projectively flat, where
$\tilde\na$ is defined by $D_XY = \tilde \na _X Y + \tilde h(X,Y)f$
for any hypersurface transverse to the position vector $f$.  A
partial converse is also true (\cite{nomizu-sasaki} or
\cite{eisenhart}):
\begin{prop} \label{centro-prop}
Any manifold equipped with   a projectively flat, torsion-free
connection with symmetric Ricci tensor is local the pull-back of the
centroaffine connection of a hypersurface under a diffeomorphism.
\end{prop}

An $\rp^n$ manifold $M$ admits a \emph{development pair} of a
developing map and holonomy.  Given a universal cover $\tilde M$ of
$M$ and a fundamental group $\pi_1M$ corresponding to a base point
$x_0$, there is a pair $({\rm dev}, {\rm hol})$ of ${\rm dev}:
\tilde M\to \rp^n$ and ${\rm hol}: \to \pgl{n+1}$ satisfying ${\rm
dev} \circ {\rm hol}(\gamma) = \gamma \circ {\rm dev}$ for all
$\gamma\in \pi_1M$. The pair ${\rm dev}$ is unique up to composition
in dev and conjugation in hol.  An analog of this theorem is valid
for all $(X,G)$ manifolds, and is due to Ehresmann.  See e.g.\
Goldman \cite{goldman-lec-notes}.  The developing map is constructed
by starting at the base point $x_0\in M$ with a projective
coordinate chart around $x_0$.  Then along any path in $\tilde M$
from a lift $\tilde x_0$ of $x_0$, there is a unique choice of
coordinate charts.  This defines the developing map ${\rm dev} :
\tilde M \to \rp^n$, and the holonomy ${\rm hol}(\gamma)$ is given
by the coordinate transformation in $\pgl{n+1}$ given by developing
the path from $x_0$ to $\gamma(x_0)$.

The developing map of an $\rp^n$ manifold closely corresponds to the
centroaffine hypersurface picture above.  If $\na$ is a projectively
flat connection satisfying the conditions of Proposition
\ref{centro-prop}, then the projection of the centroaffine
hypersurface above from $\re^{n+1}\setminus \{0\}$ to $\rp^n$
coincides with the developing map.

To each properly convex $\rp^n$ manifold, there is a dual manifold
modeled on the dual projective space $\rp_n$.  In particular, if $M
= \Omega / \Gamma$, where $\Omega$ is a bounded domain in
$\re^n\subset \rp^n$ corresponding to the regular convex cone
$\mathcal C\subset \re^{n+1}$.  Then $\Omega^*$ can be taken as the
projection to $\rp_n$ of the dual cone $\mathcal C^*\subset \re_n$.

By the uniqueness and invariance of hyperbolic affine spheres,
together with the duality result, we have

\begin{prop}
The conormal map between dual hyperbolic affine spheres provides a
natural map from any properly convex $\rp^n$ manifold $M$ to its
projective dual manifold $M^*$.  This map is an isometry with
respect to the affine metrics and it takes the Blaschke connection
on $M$ to the conjugate of the Blaschke connection on $M^*$.
\end{prop}

If a properly convex $\rp^n$ manifold $M$ is written as $M= \Omega /
\rho(\pi_1M)$ for $\Omega\subset\re^n$ a bounded convex domain and
$\rho\!:\pi_1M \to \pgl{n+1}$ a holonomy representation of the
fundamental group, then the \emph{dual $\rp^n$ structure} on $M$ is
given by $\Omega^* / \chi(\rho(\pi_1M))$ for $\chi$ the map from the
projective linear group on $\re^{n+1}$ to the projective group on
its dual $\re_{n+1}$ given $\chi (\gamma) = (\gamma^\top)^{-1}$. The
previous proposition then shows that the conormal map of hyperbolic
affine spheres provides a natural identification of $M$ with its
dual $\rp^n$ manifold.

Given a convex bounded domain $\Omega\subset \re^n\subset\rp^n$
corresponds to a convex cone $\mathcal C\subset \re^{n+1}$.  For the
unique hyperbolic affine sphere $H$ asymptotic to the boundary of
$\mathcal C$ then the projective quotient map
$\re^{n+1}\setminus\{0\} \to \rp^n$ induces a diffeomorphism
$H\to\Omega$. The affine invariants of $H$ then descend to any
projective quotient of $\Omega$, and so provide invariants on any
properly convex $\rp^n$ manifold.  See e.g.\ \cite{loftin01}
\begin{prop}
On any properly convex $\rp^n$ manifold $M$, the unique hyperbolic
affine sphere asymptotic to the cone over the universal cover of $M$
determines the following data:  a complete Riemannian metric (the
affine metric), and two canonical projectively flat connections (the
affine connection $\na$ and the conjugate connection
$\na^*=\na-2C$), representing the $\rp^n$ structure and the dual
$\rp^n$ structure.
\end{prop}
\begin{proof}
The duality theorem above \ref{dual-conj} shows that dual hyperbolic
affine sphere has affine connection equal to the conjugate
connection $\na^*=\na-2C$ for $C$ the cubic form.  Both $\na$ and
$\na^*$ are projectively flat, since they are centroaffine
connections on hypersurfaces in $\re^{n+1}$.
\end{proof}

Jaejeong Lee has recently proved the following theorem by using
hyperbolic affine spheres and their duality structure under the
conormal map.
\begin{thm}
Every compact properly convex $\rp^n$ manifold $M$ has a fundamental
domain given by a bounded polytope in $\re^n\subset\rp^n$.
\end{thm}

The case of convex $\rp^2$ structures is particularly rich, since
\c{T}i\c{t}eica's equation provides a link between hyperbolic affine
spheres and holomorphic cubic differentials on surfaces.

As in Section \ref{2-dim-sec} above, Wang \cite{wang91} provides a
version of \c{T}i\c{t}eica's developing map for two-dimensional
hyperbolic affine spheres.  Given a Riemann surface $\Sigma$ with a
holomorphic cubic differential $U$ and a complete conformal metric
$e^uh$ satisfying
$$\Delta u + 4\|U\|^2 e^{-2u} - 2 e^u - 2\kappa$$
as in Section \ref{2-dim-sec} above, the affine structure equations
produce a map from the universal cover $\tilde \Sigma \to \re^3$. On
a compact Riemann surface equipped with a hyperbolic background
metric, Wang's equation has a unique solution for each pair
$(\Sigma,U)$ \cite{wang91,labourie97,loftin01}.

Thus we have the following theorem, due independently to  Labourie
\cite{labourie97,labourie07} and the author \cite{loftin01}:

\begin{thm}
 \label{rp2-hol-cubic}
On a closed oriented surface of genus $g\ge2$, a convex $\rp^2$
structure is equivalent to a pair $(\Sigma,U)$ consisting of a
conformal structure $\Sigma$ on the surface and a holomorphic
cubic differential on $\Sigma$.
\end{thm}

One can define the deformation space of convex $\rp^2$
structures by analogously to Teichm\"uller space: On a closed
oriented surface $R$ of genus $g\ge2$, the deformation space
$\mathcal G_g$ of convex $\rp^2$ structures can be defined as
the set of equivalence classes of pairs $[S,f]$, where $S$ is a
convex $\rp^2$ surface, $f\!: R \to S$ is an
orientation-preserving diffeomorphism, and $(S,f) \sim (S',f')$
if and only if $f'\circ f^{-1}$ is homotopic to a
diffeomorphism from $S\to S'$ which preserves the projective
structure.  See e.g.\ Goldman \cite{goldman90a}.

The Riemann-Roch theorem then shows that the deformation space of
convex $\rp^2$ structures has the structure of the total space of
the holomorphic vector bundle over the Teichm\"uller space of
conformal structures on the surface whose fibers are the vector
space of holomorphic cubic differentials.

\begin{cor}[Goldman \cite{goldman90a}] The deformation space of
convex $\rp^2$ structures on a closed oriented surface of genus
$g\ge 2$ is a cell of real dimension $16g-16$.
\end{cor}

We may define the moduli space of convex $\rp^2$ structures as
the space of all projective diffeomorphism classes of convex
$\rp^2$ structures on a closed oriented surface $R$ of genus
$g\ge2$.  As in the case with Teichm\"uller space, this moduli
space of convex $\rp^2$ structures is the quotient of the
deformation space $\mathcal G_g$ by the mapping class group.
Theorem \ref{rp2-hol-cubic} allows us to describe this space as
an (orbifold) vector bundle over the moduli space of Riemann
surfaces of genus $g$.

In general, it is not easy to determine explicitly the $\rp^2$
holonomy determined by a pair $(\Sigma,U)$ for $\Sigma$ a complex
structure and $U$ a cubic differential.  But in some limiting cases,
we can determine some information.

We may determine the $\rp^2$ holonomy in  limits of pairs
$(\Sigma_t,U_t)$ for $\Sigma_t\to\Sigma_\infty$ a point in the
Deligne-Mumford compactification of the moduli space of Riemann
surfaces.  At such a limit $\Sigma_\infty$, one or more necks of the
Riemann surface are pinched to nodes.  Holomorphic cubic
differentials on $\Sigma_t$ then naturally have limits as regular
cubic differentials on $\Sigma_\infty$, which are allowed to have
poles of order 3 at each puncture.  We may define the residue of the
cubic differential to be the $dz^3/z^3$ coefficient of a regular
cubic differential with pole at $z=0$.  This residue is invariant
under the choice of local holomorphic coordinate $z$.    For regular
cubic differentials, the residues across each puncture must sum to
zero. In \cite{loftin02c}, around sufficiently small loops around
each puncture, the holonomy determined by the structure equations
(\ref{gauss-eq}-\ref{weingarten-eq}) is an ODE which approaches
$$\partial_x X = A X,$$ for a frame $X$ and a constant matrix $A$.
This determines the conjugacy class as that of $\exp (A)$, at least
when the holonomy has distinct eigenvalues.  This gives an explicit
relationship between the cubic differential and Goldman's analog of
Fenchel-Nielsen's length coordinates \cite{goldman90a}.

An $\rp^2$ structure on a surface determines (up to conjugacy) a
holonomy representation from $\pi_1 R \to \Sl3$. Using the theory of
Higgs bundles, Hitchin has identified a component of the
representation space of the fundamental group of a closed surface
$R$ of genus $g\ge2$ into $\Sl3$ (and into any split real form of a
semisimple Lie group) \cite{hitchin92}. Fixing a Riemann surface
structure $\Sigma$ on $R$, Hitchin identifies the component of the
representation space with $H^0 (\Sigma,K^2) \oplus H^0(\Sigma, K^3)$
the space of pairs of quadratic and cubic differentials over
$\Sigma$.  Choi-Goldman \cite{choi-goldman93} show that the holonomy
map identifies the space of convex $\rp^2$ structures on $R$ with
Hitchin's component of the representation space.  Under this map
Labourie \cite{labourie07} has identified Hitchin's cubic
differential with (a constant multiple of) C.P.\ Wang's cubic
differential $U$ coming from the hyperbolic affine sphere:  The
$\rp^2$ structure coming from the affine sphere determined by
$(\Sigma,-\frac1 {12} U)$ has holonomy determined by Hitchin's Higgs
bundle representation from the Riemann surface $\Sigma$ with
quadratic differential 0 and cubic differential $U$.

Labourie \cite{labourie07} also notes that a hyperbolic affine
sphere $L\subset\re^3$ naturally gives rise to a harmonic map into
the symmetric space of metrics on $\re^3$, which may be identified
with $\Sl3/\so3$.  If $L\subset\re^3$ is a hyperbolic affine sphere
with center 0 and affine mean curvature $-1$, and $p\in L$, then we
may define the \emph{Blaschke lift} $$\mathcal B\!: L \to
\Sl3/\so3$$ as the metric on $\re^3$ given by the orthogonal direct
sum of the affine metric on the tangent space $T_pL$ and the metric
on the line $\langle p \rangle$ for which the position vector $p$
has norm 1. Then with respect to the affine metric on $L$, $\mathcal
B$ is a harmonic map.  Under a projective quotient of $L$, $\mathcal
B$ is a twisted harmonic map (a section of a bundle).  Such a
harmonic metric is a natural foundation of the theory of Higgs
bundles (see Corlette \cite{corlette88}), and provides a direct link
between affine spheres and Hitchin's Higgs bundle theory.

(We also note that if the affine sphere $L$ is globally asymptotic
to the boundary of a convex cone $\mathcal C \subset \re^3$, then
the Blaschke lift on $L$ is the restriction of Cheng-Yau's complete
affine K\"ahler-Einstein metric on $\mathcal C$---see Sasaki
\cite{sasaki85} and the discussion in the next section.)



\section{Affine Manifolds} \label{affine-manifolds}

An affine structure on a smooth manifold $M$ is provided by a
maximal atlas of coordinate charts in $\re^n$ with locally constant
affine transition maps.  Thus, in terms of Thurston's notation, an
affine structure is the structure of an $(X,G)$ manifold with
$X=\re^n$ and
$$ G = \mathbf{Aff}(n,\re) = \{\Phi\!: \re^n \to \re^n, \Phi(x) =
Ax+b\}.$$  Equivalently, an affine structure is provided by a flat,
torsion-free connection $\na$ on the tangent bundle: In this case,
the geodesics of $\na$ are straight line segments in the affine
coordinate charts on $M$.

Affine manifolds are related to parabolic affine spheres through the
Monge-Amp\`ere equation $\det\phi_{ij} = 1$ (recall the graph of
such a $\phi$ is a parabolic affine sphere).  The Blaschke
connection of a parabolic affine sphere is torsion-free and flat.
Moreover, under an affine change of coordinates, the Hessian of a
function $\phi$ transforms as a tensor.  In particular, a convex
function $\phi$ defines a Riemannian metric
$\frac{\partial^2\phi}{\partial x^i
\partial x^j} dx^idx^j$ on an affine manifold.  A Riemannian metric
on an affine manifold which is locally of this form is called an
affine K\"ahler, or Hessian, metric.  If $M$ admits a $\na$-parallel
volume form $\nu$ (in other words, the affine holonomy is the
special affine group $\sa{n}$), then the Monge-Amp\`ere equation
$\det \phi_{ij} = \nu^2$ can be written as $\det \phi_{ij} = 1$ in
special affine coordinates.

The tangent bundle of an affine manifold naturally carries a complex
structure.  Let $x\mapsto Ax+b$ be an affine change of coordinates
for affine coordinates $x=(x^1,\dots,x^n)$, and let
$y=(y^1,\dots,y^n)$ denote frame coordinates on the tangent space by
representing tangent vectors as $y^i\frac{\partial}{\partial x^i}$.
Then $z^i = x^i + \sqrt{-1}y^i$ form complex coordinates on the
tangent bundle via the gluing map $z\mapsto Az+b$.  We may denote
the tangent bundle $TM$ with this complex structure as $M^\co$.
(Note this construction can also be seen as gluing together tube
domains of the form $\Omega + \sqrt{-1}\re^n \subset \co^n$ where
$\Omega$ is an affine coordinate chart.) An affine K\"ahler metric
on $M$ induces a K\"ahler metric on $M^\co$, and if the metric
satisfies the Monge-Amp\`ere equation $\det \phi_{ij}=1$, induces a
Ricci-flat K\"ahler metric on $M^\co$.

This picture is important to the conjecture of Strominger-Yau-Zaslow
\cite{syz} about mirror symmetry in string theory.  At least for
Calabi-Yau manifolds at or near certain degenerate limits in their
moduli space (called large-complex-structure limits), a Calabi-Yau
manifold $N$ of dimension $n$ is conjectured to be the total space
of a singular fibration $\pi\!:N^\circ\to M$, where $M$ is
conjectured to consist an affine manifold $M^{\rm reg}$ of dimension
$n$, together with a singular locus of real codimension two.  Over
the regular points $M^{\rm reg}$, the fibers of $\pi$ are
conjectured to be special Lagrangian tori in $N$. $N^{\rm reg } =
\pi^{-1} M^{\rm reg}$ admits a Calabi-Yau structure as a quotient of
the tangent bundle structure above: $M^{\rm reg}$ is an affine
manifold equipped with an affine flat connection $\na$ and an affine
K\"ahler metric satisfying the Monge-Amp\`ere equation $\det u_{ij}
= 1$. Therefore, the tangent bundle $T M^{\rm reg}$ admits a
Calabi-Yau metric.  If the linear part of the affine holonomy of
$M^{\rm reg}$ is integral---conjugate to a map in
$\mathbf{SL}(n,\mathbb Z)$---then we may take the quotient $TM^{\rm
reg}/ \Lambda = N^{\rm reg}$ for $\Lambda$ a $\na$-invariant lattice
in $TM^{\rm reg}$.  This provides the Calabi-Yau structure on
$N^{\rm reg}$.

The mirror conjecture of Strominger-Yau-Zaslow \cite{syz}, at
least in a simplified form, is that the mirror Calabi-Yau can
be determined by the Legendre transform of the local semi-flat
Calabi Yau potential $u$ on the affine manifold $M^{\rm reg}$.
This is the same as the duality between parabolic affine
spheres we outline in Propositions \ref{dual-affine-sphere} and
\ref{dual-conj} above: the semi-flat Calabi-Yau metrics on the
$M^{\rm reg}$ and its dual $M^{{\rm reg},*}$ are isometric,
while the affine flat connections are conjugate to each other
with respect to the affine metric.  The conjecture of
\cite{syz} also requires the torus quotients to be replaced by
their dual tori, and there should also be instanton correction
terms, which we will not discuss in this paper.

The structure of the singular locus $S = M\setminus M^{\rm reg}$ is
more complicated, but the picture is largely understood in dimension
two due to the work of Gross-Wilson \cite{gross-wilson00}.  In this
case, K3 surfaces, when properly rescaled, degenerate to a semi-flat
Calabi-Yau structure on $S^2$ minus 24 points.  The author has
produced many semi-flat Calabi-Yau structures on $\cp^1$ minus a
finite number of points by using solving the PDE of \c{T}i\c{t}eica
and Simon-Wang  $$\Delta u + 4 e^{-2u} \|U\|^2 - 2\kappa = 0$$ for
$U$ a holomorphic cubic differential on $\cp^1$ with poles of order
one \cite{loftin04}. The affine metric and holonomy are
asymptotically the same as those studied in \cite{gross-wilson00}.
Again, it is more difficult from this point of view to determine
from the cubic differential the full affine holonomy from $\pi_1
M^{\rm reg} \to \sa2$.

In dimension three, the base manifold $M^{\rm reg}$ is
conjectured to be a three-manifold minus a graph.  Generically,
one can assume the graph to have trivalent vertices, and a
fundamental problem is to produce nontrivial semi-flat
Calabi-Yau structures locally near a trivalent vertex (on a
ball in $\re^3$ minus the ``Y" vertex of a  graph).  In
\cite{lyz05,lyz-erratum}, we construct such metrics by assuming
the potential is homogeneous in a radial direction, thereby
reducing the problem to an equation on the surface $S^2$ minus
3 points.  We give two constructions in
\cite{lyz05,lyz-erratum}: For any cubic differential $U$ with
three poles of order $\le3$ on $\cp^1$, a version of
\c{T}i\c{t}eica's equation (for an appropriate background
metric)
$$ \Delta \eta + 4\|U\|^2 e^{-2\eta} - 2e^\eta - 2\kappa = 0$$
can be solved to produce a hyperbolic affine sphere  structure
on $\cp^1$ minus the pole set \cite{loftin02c}. Also, for a
cubic differential $U$ with three poles of order 2 on $\cp^1$,
we can solve the corresponding equation for an elliptic affine
sphere
$$ \Delta \eta + 4\|U\|^2 e^{-2\eta} + 2e^\eta - 2\kappa = 0$$
as long as $U$ is nonzero and small \cite{lyz-erratum}.  Then a
result of Baues-Cort\'es \cite{baues-cortes03} produces a semi-flat
Calabi-Yau structure on a ball in $\re^3$ minus the ``Y" vertex of a
graph.  There is also a construction of Zharkov \cite{zharkov04} in
which the holonomy is determined and an affine K\"ahler metric
produced, but the Monge-Amp\`ere equation is not satisfied.

Computing the affine holonomy of such solutions is still open: This
amounts to computing the projective holonomy of the solutions on
$S^2$ minus the poles of $U$.  The conjugacy class of the holonomy
around free loops around each puncture can generally be calculated,
but the problem for more complicated paths seems much harder.  See
the discussion above in Section \ref{2-dim-sec} on surfaces and
\ref{proj-man-sec} on projective structures.  Perhaps a place to
start in terms of more global holonomy calculations is the
following: Consider $U$ a meromorphic cubic differential with three
poles of order 2 on $\cp^1$.  Then, by an observation of Robert
Bryant, the parabolic affine sphere metric is given by $$\frac
{|U|^2} {m^2},$$ where $m$ is the complete hyperbolic metric on
$\cp^1$ minus the pole set of $U$. Since both the metric and cubic
differential can be made reasonably explicit in this case, it is
perhaps tractable to find the affine holonomy and developing map for
the affine structure on $\cp^1$ minus 3 points determined by this
parabolic affine sphere structure.  The Weierstrass formula for
parabolic affine spheres should help.

Compact K\"ahler affine manifolds were studied by Cheng-Yau in
\cite{cheng-yau82}, and also by Shima---see e.g.\ \cite{shima78} and
\cite{shima07}.  Cheng-Yau's work on the analogs of
K\"ahler-Einstein metrics on affine manifolds are strongly related
to affine spheres, and so we focus on \cite{cheng-yau82}.

We have the following theorem of Cheng-Yau \cite{cheng-yau82}
\begin{thm}
Let $M$ be a compact K\"ahler affine manifold which admits a
covariant constant volume form $\nu$.  Then for every smooth affine
K\"ahler metric $g$, there is a positive constant $c$ and a smooth
function $u$ so that  $$\det \left( g_{ij} + \frac{\partial ^2u}{
\partial x^i \partial x^j} \right) = c\, \nu ^2, \qquad g_{ij} + \frac{\partial ^2u}{
\partial x^i \partial x^j}>0.$$  In other words, the tensor $g+ \na
du$ is a Riemannian metric whose volume form is $\sqrt c \,\nu$.
\end{thm}
Note this theorem produces a Calabi-Yau metric on the tangent bundle
$M^\co$, and the proof uses  in an essential way Yau's estimates for
producing K\"ahler-Einstein metrics \cite{yau78}.
\begin{rem}
Cheng-Yau actually prove a more general result that given any volume
form $V$ on a compact special affine K\"ahler manifold, there are a
constant $c$ and a function $u$ so that $g+\na du>0$ and  $\det
(g_{ij} + u_{ij} ) =c\,V^2$.  Later, Delano\"e proved an analogous
theorem on any compact affine K\"ahler manifold which does not
necessarily admit a parallel volume form \cite{delanoe89}.
\end{rem}

\begin{cor} (Cheng-Yau \cite{cheng-yau82})
Every compact affine K\"ahler manifold which admits an invariant
volume form also admits a flat affine K\"ahler metric.
\end{cor}

The corollary follows from the theorem of Calabi \cite{calabi72}
that the cubic form on an affine-complete parabolic affine sphere
must vanish.  This implies the affine metric is flat.

For any volume form $V$ on an affine manifold $M$, $\na d\log V$ is
a symmetric $(0,2)$ tensor.  The analogous statement for complex
manifolds is that $\partial\bar\partial \log $ of any volume form is
a $(1,1)$ form.  Cheng-Yau \cite{cheng-yau82} prove the following
analog of Aubin and Yau's theorem on K\"ahler-Einstein metrics with
negative Ricci curvature:
\begin{thm}
If $M$ is a compact affine manifold so that $\na d\log V>0$, then
there exists a volume form $\tilde V$ on $M$ so that $$\det \left(
\frac{\partial^2 \log \tilde V} {\partial x^i \partial x^j} \right)
= \tilde V^2.$$
\end{thm}
The resulting metric $\na d\log V$ is the restriction of a
K\"ahler-Einstein metric of negative Ricci curvature on $M^\co$, and
thus may be called affine K\"ahler-Einstein metrics or Cheng-Yau
metrics.

In this case, results of Koszul \cite{koszul65}  and Vey
\cite{vey70} show that $M$ is an affine quotient of a convex cone
containing no lines (since the closed one-form $\alpha=d\log V$
satisfies $\na\alpha>0$).  Cheng-Yau also prove in
\cite{cheng-yau82}

\begin{thm}
On each convex cone $\mathcal C \subset \re^{n+1}$ which contains no
lines. Then there is a unique solution to $$ \det u_{ij} = e^{2u},
\quad u = \infty \mbox{ on } \partial \mathcal C, \quad [u_{ij}] >
0.$$  The resulting metric $u_{ij}\,dx^idx^j$ is the complete affine
K\"ahler metric generating the complete K\"ahler-Einstein metric on
the tube domain $\mathcal C + \sqrt{-1} \re^{n+1}$.
\end{thm}

 These results of Cheng-Yau are also related to affine spheres due to a result of
Sasaki \cite{sasaki85} (see also \cite{loftin02b}).  First of all,
the result of Koszul and the invariance of K\"ahler-Einstein metrics
shows that there are affine coordinates on $M$ so that $\tilde V$ is
homogeneous of degree $-n-1$.  In these coordinates, Sasaki shows
the level sets of $\tilde V$ are hyperbolic affine spheres.

The result of \cite{sasaki85} shows that two of the invariant
structures on convex cones determined by Cheng-Yau are equivalent:
In terms of the natural affine coordinates on $\mathcal C$ for which
the vertex of the cone is the origin, the level sets of the volume
of Cheng-Yau's affine K\"ahler-Einstein metric are hyperbolic affine
spheres.  Moreover, under the foliation of the cone by homothetic
copies of the hyperbolic affine sphere, the affine K\"ahler-Einstein
metric is a metric product of the affine metric on the hyperbolic
affine sphere and the flat metric $dr^2/r^2$ on the radial parameter
in $\re^+$.  Cheng-Yau's affine K\"ahler-Einstein metric is
invariant under any linear automorphism of the cone, and so it
descends to any affine quotient manifold.

There are other invariant affine K\"ahler metrics on regular convex
cones. The \emph{characteristic function} on $\mathcal C\subset
\re^n$ is defined by
$$\psi(x) = \int_{\mathcal C^*} e^{-\langle x,x^*\rangle} dx^*$$ for
$\mathcal C^*$ the dual cone of $\mathcal C$.  Then Koszul-Vinberg's
metric $\na d\log\psi$ is an affine K\"ahler metric on $\Omega$
invariant under linear automorphisms of $\mathcal C$. Thus this
metric too descends to any affine manifold which is a quotient of
$\mathcal C$ by a group of linear automorphisms acting discretely
and properly discontinuously. One can think of Koszul-Vinberg's
metric as an analog of the Bergman metric in complex geometry, but
this metric is not, for general non-homogeneous cones, the
restriction of the Bergman metric on the tube domain $\Omega +
\sqrt{-1} \mathbb R^n$ (see e.g.\ \cite{sasaki85} for an
investigation of these metrics).

The level sets of the volume form of these affine K\"ahler metrics
are also invariant hypersurfaces asymptotic to the boundary of the
convex cone $\mathcal C$, though they do not in general seem to have
the same duality property under the conormal map as the hyperbolic
affine sphere. See e.g.\ Darvishzadeh-Goldman \cite{d-goldman} for
applications to convex $\rp^2$ surfaces.

\section{Affine maximal hypersurfaces}

A natural generalization of the concept of a parabolic affine sphere
is an affine maximal hypersurface. The condition for a parabolic
affine sphere is that the affine shape operator $S$ vanishes, while
a hypersurface is affine maximal if the affine mean curvature
$H=\frac1n\mbox{tr}\,S = 0$.

The volume induced by affine metric provides a natural
equiaffine-invariant functional on all smooth strictly convex
hypersurfaces.  The first variation of this functional is
$$H=0$$ for $H = \frac1n\mbox{tr} S$ the affine mean curvature.  By analogy
with the Riemannian case, Blaschke then referred to these
stationary hypersurfaces as affine minimal \cite{blaschke}.
Much later, Calabi found that the second variation is negative
in many cases. Therefore, such a stationary hypersurface is a
local maximum for the volume functional, and thus we now call
stationary \emph{affine maximal hypersurfaces}.  Calabi proved
\cite{calabi82}
\begin{thm}
Let $L$ be an affine maximal hypersurface in $\re^{n+1}$. Then the
second variation under any compactly supported interior deformation
is negative if either
\begin{itemize}
\item $n=2$.
\item For any $n$, if $L$ is locally a graph of the form $x^{n+1} =
f(x^1,\dots,x^n)$. Moreover, $L$ is a global maximum among such
variations in this case.
\end{itemize}
\end{thm}

For $\Omega\subset\re^n$ a domain and a $u\!:\Omega\to\re$ a
smooth convex function, the graph $(x,u(x))$ is an affine
maximal hypersurface if and only if the fourth-order PDE
$$U^{ij} D_{ij} [(\det u_{ij})^{-\frac{n+1}{n+2}} ] = 0$$
is satisfied, where $[U^{ij}]$ is the cofactor matrix of the
Hessian matrix $[u_{ij}]$.

Chern conjectured that any properly embedded affine maximal surface
in $\re^3$ must be an elliptic paraboloid \cite{chern79}.  This
extension of J\"orgens's Theorem was proved by Trudinger-Wang
\cite{trudinger-wang00}. The corresponding problem is open in higher
dimension, but non-smooth viscosity solutions in dimension $n\ge10$
are presented in \cite{trudinger-wang00}.

\begin{thm}\cite{trudinger-wang00}
For $\Omega$ a domain in $\re^2$, and $u\!:\Omega\to\re$ is a smooth
convex function whose graph $(x,u(x))$ is a properly embedded affine
maximal surface in $\re^3$, then $u$ is a quadratic function.
\end{thm}

In order to prove this theorem, Trudinger-Wang use the
estimates of Caffarelli-Guti\'errez
\cite{caffarelli-gutierrez97} on solutions of the linearized
Monge-Amp\`ere equation.  Earlier, Calabi had proved that any
affine maximal surface in $\re^3$ which is both Euclidean and
affine complete must be an elliptic paraboloid \cite{calabi88}.
In fact, affine completeness implies Euclidean completeness for
hypersurfaces in $\re^{n+1}$ for $n\ge2$ by
\cite{trudinger-wang02}, and so we have the following result,
also proved by Jia-Li \cite{jia-li01}:

\begin{thm}
Any affine-complete maximal surface in $\re^3$ is an elliptic
paraboloid.
\end{thm}

We also remark that affine maximal surfaces are important
examples in the theory of integrable systems, with
Chern-Terng's construction of B\"acklund transformations for
them \cite{chern-terng80}.  The Weierstrass formula for affine
maximal surfaces is given by Calabi \cite{calabi88}, Terng
\cite{terng83}, and Li \cite{li89}.  Each affine maximal
surface is given by a holomorphic curve $Z\to \co^3$, and the
parametrization may be recovered by
$$ -\sqrt{-1}\left( Z \times \bar Z + \int Z \times dZ - \int
\bar Z \times d\bar Z\right)$$ for $\times$ the cross product.

\section{Affine normal flow}
Affine spheres are the solitons (self-similar solutions) to the
affine normal flow.   If a hypersurface $L$ is parametrized
locally by an immersion $f$, the affine normal flow is
$$\frac {\partial f }{\partial t} = \xi$$ for $\xi$ the affine
normal.  Hyperbolic, parabolic, and elliptic affine spheres are then
respectively expanding, translating, and contracting soliton
solutions for the affine normal flow.  Even though the affine normal
$\xi$ is a third-order invariant of $f$, the affine normal flow is
equivalent to the second-order parabolic flow of the $\frac1{n+2}$
power of the Gauss curvature (since $\xi = K^{\frac1{n+2}} \nu$ plus
a tangential part for $\nu$ the Euclidean unit normal).

Chow \cite{chow85} shows that compact smooth strictly convex
initial hypersurfaces in $\re^{n+1}$ converge in finite time
under the affine normal flow, and Andrews \cite{andrews96}
proves that the rescaled limit of such contracting solutions is
an ellipsoid.  Andrews \cite{andrews00} also shows that
arbitrary compact convex initial hypersurfaces are
instantaneously regularized under the affine normal flow.

Recently, Tsui and the author \cite{loftin-tsui08} extended the
affine normal flow to noncompact convex initial hypersurfaces.
One of the consequences is a parabolic proof of Cheng-Yau's
theorem on the existence of hyperbolic affine spheres:
\begin{thm}
Given any open convex cone $\mathcal C \subset \re^{n+1}$ which
contains no lines, the affine normal flow evolves the initial
hypersurface $\partial C$ to homothetically expanding copies of
the hyperbolic affine sphere asymptotic to $\partial C$.
\end{thm}
This theorem recovers Cheng-Yau's solution to the
Monge-Amp\`ere equation Dirichlet problem
$$\det u_{ij} = \left(-\frac1u\right)^{n+2},\qquad
u|_{\partial \Omega} = 0$$ for $\Omega\subset\re^n$ any convex
bounded domain\cite{cheng-yau77}, together with the result of
Gigena \cite{gigena81} and Sasaki \cite{sasaki80} that the
solution to such an equation is a hyperbolic sphere asymptotic
to the cone over $\Omega$.  In
\cite{loftin-tsui08,loftin-tsui07}, we also classify all
ancient solutions to the affine normal flow, showing them to be
either ellipsoids or paraboloids.

\bibliographystyle{abbrv}
\bibliography{thesis}

\end{document}